\documentclass{amsart}

\usepackage{amsmath}
\usepackage{amssymb}

\newtheorem{theorem}{Theorem}[section]
\newtheorem{lemma}[theorem]{Lemma}
\newtheorem{remark}[theorem]{Remark}

\numberwithin{equation}{section}

\begin{document}

\title{Pseudoconvexity is a two-dimensional phenomenon}
\author{Robert Jacobson}
\address{Department of Mathematics, Texas A\&M University, College Station, TX 77843-3368}
\email{rjacobsn@math.tamu.edu}


\date{July 2, 2009}

\begin{abstract}
An open set in $\mathbb{C}^n$ is pseudoconvex if and only if its intersection with every affine subspace of complex dimension two as seen as an open set in $\mathbb{C}^2$ is pseudoconvex.
\end{abstract}

\maketitle

\section{Introduction}
It is well known that complex analysis in several variables differs significantly from complex analysis in a single variable. The zeroes of a nonconstant holomorphic function of a single variable are always isolated, whereas the zeroes of a holomorphic function in several variables are never isolated. In one dimension the Riemann Mapping Theorem says that every simply connected bounded open set is biholomorphic to the unit ball, whereas in several dimensions this is not the case. In one dimension any open set supports a holomorphic function which cannot be holomorphically extended beyond any boundary point, whereas this is not always the case in several dimensions. An open set that does support such a function is \emph{pseudoconvex}.

In several complex variables some phenomena can be understood using one-dimensional techniques. For example, a function is holomorphic in $\mathbb{C}^n$ if it is holomorphic on every one-dimensional affine subspace, that is, on every complex line, by Hartogs's Theorem. Properties which cannot be detected using only one-dimensional data are in a sense more central to the study of several complex variables. For example, the notion of pseudoconvexity plays a central role in complex analysis in several variables. Pseudoconvexity cannot be detected by taking one-dimensional slices of an open set in $\mathbb{C}^n$ for $n\geq 2$, because, as noted above, every one-dimensional slice is pseudoconvex in $\mathbb{C}$. 

But how dependent is the notion of pseudoconvexity on dimension? Is it true that an open set in $\mathbb{C}^n$ is pseudoconvex if and only if its intersection with every two-dimensional affine subspace is pseudoconvex? Surprisingly, the affirmative answer to this question seems to be absent from the literature on the subject.

\begin{theorem}\label{maintheorem}
An open set $\Omega \subset \mathbb{C}^n$ is pseudoconvex if and only if the intersection of $\Omega$ with every affine subspace of complex dimension two as seen as a set in $\mathbb{C}^2$ is pseudoconvex.
\end{theorem}

\section{Background and Notation}

We begin with a few classical characterizations of pseudoconvex open sets. (For more on characterizing pseudoconvex sets see~\cite{krantz}.) Let $\Omega$ be an open set in $\mathbb{C}^n$. The following are equivalent:
\begin{enumerate}
\item There exists a holomorphic function on $\Omega$ which cannot be extended holomorphically past any boundary point of $\Omega$. A domain (i.e., a connected open set) satisfying this property is often called a domain of holomorphy. 
\item For any family of holomorphic maps of the unit disc in $\mathbb{C}$ into $\Omega$, if the images of the boundaries of these maps are all contained in a single compact subset of $\Omega$, then the whole of the images of the maps are also contained in a single compact subset of $\Omega$. This property is often referred to by the German word \emph{Kontinuit\"atssatz}, meaning literally continuity theorem.
\end{enumerate}
We say that $\Omega$ has $C^2$ boundary if there exists a  $C^2$ \emph{defining function} for $\Omega$, that is, a twice continuously differentiable function $\rho \colon \mathbb{C}^n \to \mathbb{R}$ such that $\Omega = \{z\in \mathbb{C}^n \mid \rho(z) < 0\}$, $\partial \Omega = \{z\in \mathbb{C}^n \mid \rho(z)=0\}$, and $\nabla \rho(M) \neq 0$ for all $M \in \partial \Omega$. In this case, the following is also equivalent to the above.
\begin{enumerate}
\setcounter{enumi}{2}
\item If $\rho \colon \mathbb{C}^n \to \mathbb{R}$ is a $C^2$ defining function for $\Omega$ then for every point $M$ in the boundary of $\Omega$,
\[
\sum_{j, k = 1}^n \frac{\partial^2 \rho}{\partial z_j \partial \overline{z}_k}(M) Z_j \overline{Z}_k \geq 0
\]
for every $Z\in \mathbb{C}^n$ satisfying
\[
\sum_{j= 1}^n \frac{\partial \rho}{\partial z_j }(M) Z_j = 0.
\]
This property is called Levi pseudoconvexity after the Italian mathematician Eugenio Elia Levi. The left-hand expression in the first displayed inequality is known as the Levi form, while vectors $Z\in \mathbb{C}^n$ satisfying the second displayed equation are called complex tangent vectors at the point~$M$. We denote by $\mathcal{T}_M(\partial \Omega)$ the set of all complex tangent vectors at~$M$. This property says that the Levi form is nonnegative on the boundary of~$\Omega$.
\end{enumerate}

From the point of view of Property~1, Theorem~\ref{maintheorem} says that if for each two-dimensional slice of $\Omega$ there exists a holomorphic function on the slice which cannot be extended beyond any boundary point of the slice, then there exists a holomorphic function on $\Omega$ which cannot be extended past any boundary point of $\Omega$. From this perspective it is not at all clear why such a theorem should exist.

Let us establish some notation. If $a, b, c\in \mathbb{C}^n$ with $b$ and $c$ linearly independent over $\mathbb{C}$, then denote by $h=h(a,b,c):=\{a+bw_1+cw_2 \mid (w_1, w_2)\in\mathbb{C}^2\}$ the complex affine subspace of complex dimension two determined by $a$, $b$, and $c$. Define $\phi \colon \mathbb{C}^2 \to \mathbb{C}^n$ by $\phi(w_1, w_2)=a +b w_1 +c  w_2$ for all $(w_1, w_2)\in \mathbb{C}^2$ (where $a, b, c$ define $h$ and are clear from context). For any set $\Omega \subset \mathbb{C}^n$, define $\Omega_h:=\phi^{-1}(\Omega)$, the slice of $\Omega$ by $h$ viewed as a set in $\mathbb{C}^2$. For any function $f\colon\mathbb{C}^n\to\mathbb{C}$, define $f_h\colon\mathbb{C}^2\to\mathbb{C}$ by $f_h:=f\circ \phi$.

\section{Pseudoconvexity of the open set implies pseudoconvexity of every slice}

It is easy to show that every intersection with a two-dimensional affine subspace of a pseudoconvex open set is pseudoconvex. One way is to use the Kontinuit\"{a}tssatz in the following way. Let $\Omega\subset\mathbb{C}^n$ be a pseudoconvex open set and fix some $h=h(a, b, c)$ as in the previous section. We assume without loss of generality that $\Omega_h\neq \varnothing$. Let $\{d_\alpha\}_{\alpha\in A}$ be a family of analytic discs in $\Omega_h$ such that $\cup_{\alpha \in A}\partial d_\alpha$ is contained in a compact subset of $\Omega_h$. Observe that $\{\phi \circ d_\alpha\}_{\alpha\in A}$ is a family of analytic discs in $\Omega$ such that $\cup_{\alpha \in A} \partial (\phi \circ d_\alpha)$ is contained in a compact subset of $\Omega$. Since $\Omega$ is pseudoconvex, $\cup_{\alpha \in A} \phi \circ d_\alpha$ is contained in a compact subset of $\Omega$. Hence, $\cup_{\alpha \in A}d_\alpha$ is contained in a compact subset of $\Omega_h$.

The other direction in Theorem~\ref{maintheorem} takes more work.

\section{The case that $\Omega$ has $C^2$ boundary}

We begin with the special case of $\Omega$ having $C^2$ boundary.

\begin{theorem}\label{specialcase}
Let $\Omega \subset \mathbb{C}^n$ be an open set with $C^2$ boundary.  If each slice of $\Omega$ by a complex affine subspace of complex dimension two is pseudoconvex when viewed as a set in $\mathbb{C}^2$, then $\Omega$ is pseudoconvex.
\end{theorem}
To prove this we first prove a lemma which tells us when a $C^2$ defining function in $\mathbb{C}^n$ induces a $C^2$ defining function on an two-dimensional affine subspace.

\begin{lemma}\label{subspacedefiningfunction}
Let $\Omega \subset \mathbb{C}^n$ be an open set with $C^2$ boundary and let $\rho$ be a $C^2$ defining function for $\Omega$. Let $h(a, b, c)$ be a complex affine subspace such that $\Omega_h \subset \mathbb{C}^2$ is nonempty. Let $\mu \in \partial \Omega_h$. If either $b$ or $c$ is not complex tangent to $\Omega$ at $\phi(\mu)$, then $\Omega_h$ has $C^2$ boundary in a neighborhood of $\mu$ with $C^2$ defining function $\rho_h$.
\end{lemma}

\begin{proof}[Proof of Lemma \ref{subspacedefiningfunction}]
Let $\mu \in \partial\Omega_h$. Observe that $\rho_h$ is defined in a neighborhood $U$ of $\mu$. That $\rho_h$ is $C^2$ in $U$ is clear from the chain rule. Also, by construction we have $U \cap \Omega_h = \{w\in U \mid \rho_h(w)<0\}$. It remains to check that the gradient of $\rho_h$ is nonzero on $U\cap\partial\Omega_h$. What we need to show is that either the $w_1$ derivative or the $w_2$ derivative of $\rho(a +bw_1 +cw_2)$ is nonzero at $\phi(\mu)$. These derivatives are $\sum_{j=1}^n \frac{\partial \rho}{\partial z_j} b_j$ and $\sum_{j=1}^n \frac{\partial \rho}{\partial z_j} c_j$ respectively, and by hypothesis one of these derivatives is nonzero. This proves the lemma.
\end{proof}

\begin{proof}[Proof of Theorem \ref{specialcase}]

Our strategy will be to show that $\Omega$ is Levi pseudoconvex, that is, that the Levi form is nonnegative on the boundary of $\Omega$. Let $\rho$ be a $C^2$ defining function for $\Omega$. Let $M\in\partial\Omega$ and let $Z\in\mathcal{T}_M(\partial\Omega)$. Choose $p_0\in\Omega$ such that $M-p_0$ is orthogonal to $\mathcal{T}_M(\partial\Omega)$. Let $a=p_0$, $b=M-p_0$, and $c=Z$ and consider $h(a, b, c)$. Observe that $b$ and $c$ are linearly independent over $\mathbb{C}$ and that $b\not\in\mathcal{T}_M(\partial\Omega)$. Hence $\rho_h$ is locally a $C^2$ defining function for $\Omega_h$ near $\phi^{-1}(M)$ by Lemma~\ref{subspacedefiningfunction}. Let $\mu:=\phi^{-1}(M)$. Define $\zeta$ as $\zeta:=\phi^{-1}(p_0+Z)-\phi^{-1}(p_0) = (0, 1)$.

We now check that $\zeta\in\mathcal{T}_{\mu}(\partial\Omega_h)$. We compute by the chain rule:
\begin{eqnarray*}
 & & \sum_{j=1}^2 \frac{\partial \rho_h}{\partial z_j}(\mu) \zeta_j
 =  \frac{\partial \rho_h}{\partial w_2} (\mu) 
 =  \frac{\partial (\rho\circ\phi)}{\partial w_2}(\mu) \\
 & = & \sum_{k=1}^n 
\left[\frac{\partial \rho}{\partial z_k}(\phi(\mu))
\cdot \left(\frac{\partial \phi}{\partial w_2} (\mu) \right)_k
+ \frac{\partial\rho}{\partial \overline{z}_k}(\phi(\mu))
\cdot \left( \frac{\partial\overline{\phi}}{\partial w_2} (\mu) \right)_k \right] \\
 & = & \sum_{k=1}^n \frac{\partial\rho}{\partial z_k}\left(\phi (\mu) \right)Z_k 
= \sum_{k=1}^n \frac{\partial\rho}{\partial z_k}\left(M\right)Z_k = 0.
\end{eqnarray*}
Hence $\zeta\in\mathcal{T}_{\mu}(\partial\Omega_h)$.

Now, $\rho_h$ is a defining function for $\Omega_h$ near $\mu$, and by hypothesis $\Omega_h$ is pseudoconvex, so
\begin{align*}
0  & \leq \sum_{j,k=1}^2 \frac{\partial^2 \rho_h}{\partial w_j \partial \overline{w}_k}(\mu) \zeta_j \overline{\zeta}_k 
 = \frac{\partial^2 \rho_h}{\partial w_2 \partial \overline{w}_2}(\mu) \\
& = \frac{\partial^2 (\rho\circ\phi)}{\partial w_2 \partial \overline{w}_2}(\mu) 
 = \sum_{\ell, m=1}^n \frac{\partial^2 \rho}{\partial z_\ell \partial \overline{z}_m} (\phi(\mu)) c_\ell \overline{c}_m \\
& = \sum_{\ell, m=1}^n \frac{\partial^2 \rho}{\partial z_\ell \partial \overline{z}_m} (M) Z_\ell \overline{Z}_m,
\end{align*}
which proves the theorem.
\end{proof}

\begin{remark}\label{specialcaseremark}
\emph{The preceding proof shows something slightly stronger than Theorem~\ref{specialcase}: if $\Omega$ has $C^2$ boundary near a point of nonpseudoconvexity $M\in \partial\Omega$ then there is a slice $\Omega_h$ of $\Omega$ such that $h$ \emph{contains} $M$ and $\Omega_h$ is not pseudoconvex.}
\end{remark}

\section{The General Case}

\begin{theorem}\label{generalcase}
Let $\Omega \subset \mathbb{C}^n$ be an open set.  If each slice of $\Omega$ by a complex affine subspace of complex dimension two is pseudoconvex when viewed as a set in $\mathbb{C}^2$, then $\Omega$ is pseudoconvex.
\end{theorem}

Proving Theorem~\ref{generalcase} will complete our proof of Theorem~\ref{maintheorem}. 

We will prove the contrapositive. The idea of the proof is that if our open set is not pseudoconvex, we can find a point of nonpseudoconvexity on the boundary of the open set such that the point is also a point of nonpseudoconvexity of an open set with smooth boundary sitting inside the original open set. We then use our result for the case of an open set with $C^2$ boundary to show that there is a slice of this smooth open set that is not pseudoconvex. But then there is a slice of our original open set such that at a point on the boundary there is an open set with smooth boundary contained in the slice with that point as a point of nonpseudoconvexity. Thus that slice of our original open set is not pseudoconvex. To execute this strategy we will use the following theorem from~\cite[p. 240]{hormander}. 

\begin{theorem}[H\"{o}rmander]\label{hormandertheorem}
Let $\Omega$ be an open set in $\mathbb{C}^n$ which is not pseudoconvex. Then there is a point $M\in \partial \Omega$, a quadratic polynomial $q \colon \mathbb{C}^n \to \mathbb{R}$, and a neighborhood $U$ of $M$ such that $q(M)=0$, $\nabla q(M)\neq 0$, and whenever $q(z) < 0$ then $z\in \Omega$ for all $z\in U$, and
\[
\sum_{j=1}^n \frac{\partial q(M)}{\partial z_j} Z_j= 0,\qquad
\sum_{j, k = 1}^n \frac{\partial^2 q(M)}{\partial z_j \partial \overline{z}_k} Z_j\overline{Z}_k < 0
\]
for some $Z \in \mathbb{C}^n$. Conversely, $\Omega$ is not pseudoconvex if there exists such a polynomial.
\end{theorem}

\begin{proof}[Proof of Theorem \ref{generalcase}]
Assume $\Omega \subset \mathbb{C}^n$ is a nonpseudoconvex open set. Then there is a point $M\in \partial \Omega$, a quadratic polynomial $q \colon \mathbb{C}^n \to \mathbb{R}$, a neighborhood $U$ of $M$, and a $Z \in \mathbb{C}^n$ as in Theorem~\ref{hormandertheorem}. Set $V:=\{z \in U \mid q(z)<0\}$. By Remark~\ref{specialcaseremark} there exists an affine subspace $h$ containing $M$ such that $V_h$ is not pseudoconvex at $\phi^{-1}(M)$. By Lemma~\ref{subspacedefiningfunction}, $\rho_h$ is a defining function for $V_h$ near $\phi^{-1}(M)$, and $\rho_h$ is a real-valued quadratic polynomial. By the converse part of Theorem~\ref{hormandertheorem} we have that $\Omega_h$ is not pseudoconvex.
\end{proof}


\begin{thebibliography}{9}

\bibitem{krantz}
  Steven G. Krantz,
  \emph{Function Theory of Several Complex Variables, Second Edition}.
  AMS Chelsea Publishing, Providence,
  2000. 
 
\bibitem{hormander}
  Lars H\"{o}rmander,
  \emph{Notions of Convexity}.
  Birkh\"{a}user, Boston,
  2007. 

\end{thebibliography}
\end{document}